\documentclass[11pt, oneside, a4paper]{article}
\usepackage{amsfonts, epsf}
\usepackage{latexsym}
\usepackage{graphicx}
\usepackage{amssymb, amsmath}
\usepackage{xcolor}
\usepackage{hyperref}
\usepackage{authblk}
\usepackage{hyperref}
\addtocontents{11}{text}
\begin{document}
%\baselineskip=11pt
%$$$$$$$$$$$$$$$$$$$$$$$$$$$$$$$$$$
 \begin{center}
\Large \textbf{Ergodicity and Perturbation Bounds for $M_t/M_t/1$ Queue with Balking, Catastrophes, Server Failures and Repairs}
 \end{center}
\bigskip
\centerline{Alexander Zeifman$^1$, Yacov Satin$^2$, Ivan Kovalev$^3$, Sherif I. Ammar$^4$} ,\centerline{\small 1 Vologda State University,}\\
\centerline{Institute of Informatics Problems, FRC CSC RAS,}\\
\centerline{ Vologda Research Center RAS,}\\
\centerline{ Moscow Center for Fundamental and Applied Mathematics, Moscow State University, Russia}
\centerline{\small 2 Vologda State University, Russia} \centerline{\small 3 Vologda State University, Russia}
\centerline{\small 4 Department of Mathematics, Faculty of Science, Menofia University, Shebin El Kom, Egypt}
\centerline{ Department of Mathematics, College of Science, Taibah
University, Saudi Arabia}
%\def\thefootnote{\fnsymbol{*}}
%===================================================================
\begin{abstract}

In this paper, we display methods for the computation of convergence and perturbation bounds for $M_t/M_t/1$ system with balking, catastrophes, server failures and repairs. Based on the logarithmic norm of linear operators, the bounds on the rate of convergence, perturbation bounds, and the main limiting characteristics of the queue-length process are obtained. Finally, we consider the application of all obtained estimates to a specific model.
%\\\\\\ \textbf{Keywords:} Two-Processor Heterogeneous System; Catastrophes; Repairs; Bounds; Rate of convergence.
\end{abstract}
%===================================
 \section{\textit{Introduction}}
$~~~~~~~$Recently, there has been a noticeable interest from researchers to study nonstationary queueing systems because they represent the actual reality of many applications in our life.
Nevertheless, we find few works around these systems, as studying these systems needs new, unrecognized mechanisms to analyze their behavior.\\
$~~~~~~~$In the current paper, we deal with nonstationary $M/M/1$ queue with balking, catastrophes, server failures and repairs. We investigate the bounds on the rate of convergence, and the perturbation bounds for the corresponding queue-length process. Such kinds of bounds give us the possibility for finding the limiting bounds for the class of close to this queue
Markovian models. We apply the approach based on the notion of logarithmic norm of an operator function, see detailed description in our recent papers \cite{Ammar2020,Zeifman2019amc}.
The motivation of the proposed system comes from having wide and many applications and contributions in many fields, one of them is the field of communication network systems. The applicability of this model can be seen in communication network systems. If there are numerous packets lined up in the system, local packets are always accepted and remote packets are less than a threshold value packets, waiting in the node for process. Then a new arrival either decides not to join the system or departs after joining the system. If this network was infected with a virus, this could lead to the loss of some packets as a result of restarting the network again, or transferring theses packets to another network.
Also, in computer systems where there are several clients (data) lined up in the system until a certain threshold value, a new arrival may decide not to enter the system after that value. Additionally, if a virus-infected data will annihilate or transfer it to other processors. These systems can be described as queueing models with catastrophes and balking. These systems can be represented as proposed queueing model.\\
$~~~~~~~$Most of the literature on the subject of the present paper focused solely on the study of stationary behavior, although this behavior is a special case, as well as the parameters in many applications are varying with time in our daily life. An example of some works that has discussed queueing systems related to the subject of this paper, we find that in \cite{Ammar2014} the author discussed the stationary behavior of a two-processor heterogeneous system with catastrophes, server failures and repairs. Kumar et al. \cite{Krishna2007} analyzed the stationary behavior for an $M/M/1$ queue with catastrophes, server failures and repairs. In \cite{Tarabia} the author has extend work which has been done by  Kumar et al. \cite{Krishna2007} for an $M/M/1$ queue with balking, catastrophes, server failures and repairs where balking occurs if and only if the system size equals or exceeds a threshold value $k$. Also, Suranga et al. \cite{Suranga2018} considered an $M/M/1$ queue with reneging, catastrophes, server failures and repairs, they obtained the explicit expression for the stationary probabilities. While in \cite{DiCrescenzo2012} Crescenzo et al. studied the stationary behavior of a double-ended queue with catastrophes and repairs.\\
$~~~~~~~$On the other hand, we find some of the works that discussed the nonstationary behavior as in \cite{Ammar2020} Ammar et al. explored the nonstationary of a two-processor heterogeneous system with catastrophes, server failures and repairs. In \cite{DiCrescenzo2018} DiCrescenzo et al. construed the a time-non-homogeneous for double-ended queue subject to catastrophes and repairs, as this is an extension of their previous work in\cite{DiCrescenzo2012}.\\
$~~~~~~~$We note from the previous literature that no paper has discussed the behavior of the proposed model and based on this observation, in this paper we examine convergence bounds for a non-stationary behavior of the proposed system. In case of stat constant parameters, our results are consistent with those found by Tarabia \cite{Tarabia}.\\
$~~~~~~~$The paper is organized as follows. In Section 2, description of the model and basic notions  are introduced. In Section 3 and Section 4, general theorems on the rate of convergence and perturbation bounds are considered, respectively. Finally, in Section 5, a specific queueing example is studied.

%==================================================
\section{\textit{Model Description and Basic Notions}}
\numberwithin{equation}{section}
$~~~~~~~$The proposed system in the current paper is an $M/M/1$ catastrophic queue involving balking, server failure and repairs.  The arrival process of customers' is Poisson process with mean $\lambda$ arrival rate during times that the server is running. Suppose that the discipline of customers is served on first-come, first-served with the service time following an exponential appropriation with mean 1/$\mu$. When the customer arrives the system, his joining to the system depends on threshold value of $k$. If the number of customers is a fewer the threshold value of $k$, they join the queue with probability one. Also, if the number of customers is more than or equal the threshold value of $k$ they join the queue with probability $\beta$ and may balk with probability $1-\beta$. System capacity is infinite. At the point when the system is inactive or busy, catastrophes happen at the service station as indicated by Poisson process of rate $\gamma$. If a failure happens on the busy server, all the system 's customers are automatically pulverized and the server is inactivated, i.e. the server fails and needs to repair it. Failed server repair times are $i.i.d$, based on an exponential distribution with $\eta$ parameter. After the server has been repaired, the system is available to provide the service of new customers.
Let $ r(t) $ be the probability of the server at instant $t$ with $r(0) = 0 $ is under repair.\\

\smallskip

Unlike previous studies (see \cite{Krishna2007,Suranga2018,Tarabia}),
we consider in the paper the non-stationary case, that is, we suppose  that
all possible transition intensities $\lambda(t)$, $\mu(t)$, $\beta(t)$, $\gamma_i(t)$, $\eta(t)$, are
non-random functions of time, which are nonnegative and locally integrable on $[0,\infty )$.

$~~~~~~~$According to the above assumptions, the system can be described by Markov process $ {X(t), t > 0} $ where $ X(t) $ refers to the number of the system customers $t$ (queue-length process) at the time.
Denote by $ p_{n}(t) = P(X(t) = n)$, $n = 0, 1, 2, 3, \dots.$
From the previous presumptions, the resulting conduct of the state probabilities are described by forward Kolmogorov system as:
\begin{equation}
r'(t) = -\eta(t) r(t) + \sum_{i=0}^{\infty}\gamma_i(t) p_{i}(t)
\label{k1}
\end{equation}
\begin{equation}
p_{0}'(t) = \mu(t) p_{1}(t) - (\lambda(t) + \gamma_0(t))p_{0}(t) + \eta(t)
r(t)\label{k2}
\end{equation}
\begin{equation}
p_{n}'(t) = \lambda(t) p_{n-1}(t) - (\lambda(t) + \gamma_n(t) + \mu(t))p_{n}(t) +
\mu(t) p_{n+1}(t), 1 \le n \le k - 1 \label{k3}\end{equation}
\begin{equation}
p_{k}'(t) = \lambda(t) p_{k-1}(t) - (\lambda(t)\beta(t) + \gamma_k(t) +
\mu(t))p_{k}(t) + \mu(t) p_{k+1}(t), n = k \label{k4}\end{equation}
\begin{equation}
p_{n}'(t) = \lambda(t)\beta(t) p_{n-1}(t) - (\lambda(t)\beta(t) + \gamma_n(t) +
\mu(t))p_{n}(t) + \mu(t) p_{n+1}(t), n > k. \label{k5}\end{equation}
%with the initial condition $ p_{m}(0) = p_{m}$

\bigskip

$~~~~~~~$Assuming that ${\bf p}\left(t\right)$ represents the state vector of probabilities at the moment $t$, 
where ${\bf p}\left(t\right) = \left(r(t),
p_{0}\left(t\right), p_{1}\left(t\right), \dots\right)^T$. 

Put $a_{ij}\left(t\right)
=  q_{ji}\left(t\right)$ for $j\neq i$ and $a_{ii}\left(t\right) =
-\sum_{j\neq i} a_{ji}\left(t\right) = -\sum_{j\neq i}
q_{ij}\left(t\right).$

We will supposed that $|a_{ii}\left(t\right)| \le  L < \infty,$  for any $i$ and almost all $t \ge 0$.

We symbolize of the $l_1$-norm of vector by $\|\cdot\|$ , $\|{x}\|=\sum|x_i|$,
$\|B\| = \sup_j \sum_i |b_{ij}|$, if $B =(b_{ij})_{i,j=0}^{\infty}$, and set of all vectors with non-negative coordinates and unit norm from $l 1$ by $\Omega$.
We have $\|A(t)\| = 2\sup_{k}\left|a_{kk}\left(t\right)\right| \le 2L$ for almost  all  $t \ge 0$.

Therefore, in the space of sequences $l_1$, we can rewrite the forward Kolmogorov system (\ref{k1})--
(\ref{k5}) as a differential equation
\begin{equation}
\frac{d {\bf p} \left(t\right)}{dt}=A\left(t\right){\bf p} \left(t \right),
\label{s1}
\end{equation}
 where  $A\left(t\right)$ is a bounded for almost all $t \ge 0$ linear operator from $l_1$ to itself and it is the respective transposed intensity matrix is generated as : {\tiny
    \begin{equation}
    \left(
    \begin{array}{ccccccccccccccc}
    -\eta(t) & \gamma_0(t) & \gamma_1(t) & \gamma_2(t) &  \cdots & \gamma_{k-1}(t) & \gamma_k(t) & \cdots& \cdots\\
    \\
    \eta (t) & -(\lambda (t)+\gamma_0 (t))   & \mu(t) & 0 &  & \cdots & \cdots & \cdots  & \cdots \\
    \\
    0 & \lambda (t) & -(\lambda (t)+\gamma_1 (t)+\mu(t)) & \mu(t)  & 0 & \cdots & \cdots &  \cdots &  \cdots\\

    \vdots & \vdots & \vdots & \vdots & \vdots & \cdots & \vdots & \cdots &  \cdots\\
    \\
    0   & 0 & \cdots & \cdots & 0 & \lambda(t)  & -(\lambda (t) \beta (t)+\gamma_k (t)+\mu(t)) & \mu(t) &  0 & \cdots \\
    \\
    0  & 0 & \cdots  & \cdots & \cdots & 0 & \lambda(t) \beta(t)  & -(\lambda (t) \beta (t)+\gamma_{k+1} (t)+\mu(t)) & \mu(t) &  \cdots \\
    \\
    \vdots & \vdots & \vdots & \vdots & \vdots & \vdots & \vdots  &  \cdots &  \cdots\\

    \end{array}
    \right).
    \label{m1}
    \end{equation}
}

\smallskip

The mathematical expectation (the mean) of $X(t)$ is symbolized by $E(t,k) = E\left\{X(t)\left|X(0)=k\right.\right\}$ if $X(0)=k$ at the moment $t$.

\smallskip

Remember that $X\left(t\right)$ is called {\it a weakly ergodic Markov chain}, if $\lim\limits_{t
\rightarrow \infty }\left\| \mathbf{p}^1\left( t\right)
-\mathbf{p}^2\left( t\right) \right\| =0$ for any initial conditions
$\mathbf{p}^1\left( 0\right) =\mathbf{p}^1\in \Omega$, $\mathbf{p}^2\left( 0\right) =\mathbf{p}^2\in \Omega$; and it has the
limiting mean $\phi(t)$, if  $ \left|E(t;k) - \phi(t) \right| \to 0
$ as $t \to \infty$ for any $k$.

\smallskip

We use in the paper the notion of  the logarithmic norm of operator function from
$l_1$ to itself, it  is calculated by the formula
\begin{equation}
 \gamma \left( B\left(t\right) \right)_{1} = \sup_i \left(b_{ii}\left(t\right) + \sum_{j\neq i} |b_{ji}\left(t\right)|\right).
 \label{k11}
\end{equation}
\noindent Moreover, the following bound holds
\begin{equation}\|U\left(t,s\right)\| \le e^{\int_s^t
\gamma\left(B(\tau)\right)\, d\tau},     \label{k12}
\end{equation}
\noindent  where  $U\left(t,s\right)$  is the Cauchy operator of the corresponding
differential equation $\frac{d{\bf x}}{dt}=B(t){\bf
x}$.

%==================================================
\section{\textit{Bounds on the rate of convergence}}
\numberwithin{equation}{section}

As we noted earlier, our method based on the notion of logarithmic
norm and the corresponding bounds for the Cauchy operator. Moreover,
for  the considered model we can use the both approaches of
\cite{Ammar2020,Zeifman2019amc}. Describe briefly these approaches.

{\bf First approach}, see \cite{Ammar2020}. Denote by $\gamma^*(t) = \inf_{n}\gamma_n(t)$, by ${\bf
g}\left(t\right)=\left(\gamma^*\left(t\right),0,0, \dots\right)^T$.  Put 
\begin{eqnarray}
a_{ij}^*\left(t\right) = \left\{
\begin{array}{ccccccc}
a_{0j}\left(t\right) - \gamma^*\left(t\right), & \mbox { if }  i= 0, \\
a_{ij}\left(t\right), & \mbox { otherwise }.
\end{array}
\right.
\label{f2}
\end{eqnarray}
Let $A^*\left(t\right)=\left (
a_{ij}^*\left(t\right)\right)_{i,j=0}^{\infty}$. Then we can consider the equation
\begin{equation}
\frac{d {\bf p} }{dt}=A^*\left( t\right) {{\bf p} }  +{\bf g}
\left(t\right), \quad t\ge 0, \label{s2}
\end{equation}
\noindent instead of (\ref{s1}),
where
$$A^*\left(t\right)=$$
{\tiny
    \begin{equation}
    \left(
    \begin{array}{ccccccccccccccc}
    -(\eta(t) + \gamma^*(t)) & \gamma_0(t) - \gamma^*(t) & \gamma_1(t) - \gamma^*(t) &  \cdots &  \cdots &  \cdots &  \cdots &  \cdots &  \cdots &  \cdots &  \cdots &  \cdots  & \cdots \\
    \\
    \eta (t) & -(\lambda (t)+\gamma_0 (t))   & \mu(t) & 0 & 0 & \cdots & 0 & 0 & 0 & 0 & 0 & \cdots & 0  & 0 & \cdots \\
    \\
    0 & \lambda (t) & -(\lambda (t)+\gamma_1 (t)+\mu(t)) & \mu(t)  & 0 & \cdots & 0 & 0 & 0 & 0 & 0 & \cdots & 0  & 0 & \cdots\\

    \vdots & \vdots & \vdots & \vdots & \vdots & \cdots & \vdots & \vdots & \vdots & \vdots & \vdots & \cdots & \vdots & \vdots & \cdots\\
    \\

    0 & 0 & 0 & 0 & 0 & \cdots & 0 & \lambda(t)  & -(\lambda (t) \beta (t)+\gamma_k (t)+\mu(t)) & \mu(t) &  0 & \cdots & 0 & 0 & \cdots\\
    \\
    0 & 0 & 0 & 0 & 0 & \cdots & 0 & 0 & \lambda(t) \beta (t)  & -(\lambda (t) \beta (t)+\gamma_{k+1} (t)+\mu(t)) & \mu(t) &  \cdots & 0 & 0 & \cdots\\
    \\
    \vdots & \vdots & \vdots & \vdots & \vdots & \vdots & \vdots & \vdots & \vdots & \vdots & \vdots & \vdots & \vdots & \vdots & \vdots \\

    \end{array}
    \right).
    \label{m2}
    \end{equation}
}

If we denote by $U^*\left(t,s\right)$ the Cauchy operator of the corresponding homogeneous equation $
\frac{d {\bf x} }{dt}=A^*\left( t\right) {{\bf x} },$ then the equation (\ref{s2}) can be solved by the formula 
\begin{equation}
{\bf p}\left(t\right)= U^*\left(t,0\right){\bf p}\left(0\right)+\int_0^t{}U^*\left(t,\tau\right){\bf g}\left(\tau\right)\,d\tau.
\label{f3}
\end{equation}

\bigskip

Let   $1 = d_0 \le d_1
\le \dots$ be positive numbers. Denote by $ \textsf{D} = diag\left(d_0,
d_1, d_2, \dots \right)$ the corresponding diagonal matrix.
Consider the auxiliary space of sequences  $l_{1\textsf{D}} = \left\{{\bf p} / \|{\bf
p}\|_{1\textsf{D}}=\|\textsf{D}{\bf p}\|_1 < \infty\right\} $, and put
\begin{equation}
\gamma_{**}(t) =  \inf_i \left(|a_{ii}^*(t)| - \sum_{j\neq i}
\frac{d_j}{d_i}a_{ji}^*(t)\right). \label{f8}
\end{equation}

\bigskip
Then, using the arguments as in \cite{Ammar2020}, we obtain the following statements.

{\bf Theorem 1.} \label{t1} In the situation of sufficiently large catastrophe rate (i. e. the following equality holds):
\begin{equation}
\int_0^\infty \gamma^*\left(t\right)\, dt = +\infty, \label{f4}
\end{equation}
the corresponding queue-length process $X\left(t\right)$ is weakly ergodic in
the uniform operator topology. Moreover, for any different initial conditions  ${{\bf p}}^{*}\left(0\right),
{{\bf p}}^{**}\left(0\right)$  and any $t \ge 0$ we have 
\begin{eqnarray}
\left\|{{\bf p}}^{*}\left(t\right)-{{\bf
p}}^{**}\left(t\right)\right\|\le e^{-\int\limits_0^t
\gamma^*\left(\tau\right)\, d\tau}\left\|{{\bf
p}}^{*}\left(0\right)-{{\bf p}}^{**}\left(0\right)\right\| \le 2
e^{-\int\limits_0^t \gamma^*\left(\tau\right)\, d\tau}. \label{f5}
\end{eqnarray}

\smallskip

{\bf Theorem 2.} \label{t2} Let  there exist a positive sequence
$\{d_i\}$, $1 = d_0 \le d_1 \le \dots$ such that,
\begin{equation}
\int_0^\infty \gamma_{**}(t)\, dt = +\infty. \label{f11}
\end{equation}
Then  $X(t)$ is weakly ergodic and the following bound on the rate
of convergence holds:
\begin{eqnarray}
\left\|{\bf p}^{*}(t)-{\bf p}^{**}(t)\right\|_{1\textsf{D}} \le
e^{-\int\limits_0^t \gamma_{**}(\tau)\, d\tau}\left\|{\bf
p}^{*}(0)-{\bf p}^{**}(0)\right\|_{1\textsf{D}}, \label{f12}
\end{eqnarray}
\noindent for any initial conditions ${\bf p}^{*}(0), {\bf
p}^{**}(0)$ and any $t \ge 0$.

\smallskip

Let $l_{1E}=\left\{{\bf p} =(r, p_0,p_1,p_2,\ldots)\right\}$ be a
space of sequences such that  $\|{\bf p}\|_{1E}=\sum_{k \ge 0} k
|p_k| <\infty$. $ \|{\bf p}\|_{1D} = \|D {\bf p}\| = \|\left(d_0 r,
d_1 p_0, d_2 p_1, \ldots \right)^T\| = d_0 r + \sum_{k \ge 0}
d_{k+1} p_k \ge \sum_{k \ge 1} k \frac{d_{k+1}}{k} p_k. $ Put $W =
\inf_{k \ge 1}\frac{d_{k+1}}{k}$. Then $W\|{\bf p}\|_{1E} \le \|{\bf
p}\|_{1\textsf{D}}$.

\medskip

{\bf Corollary 1.}   \label{cor1} Let, under assumptions of Theorem 2, in addition inequality $W > 0$, holds. 
Then  there exists  the limiting mathematical expectation, say $\phi(t)=E(t,0)$, and the inequality 
\begin{equation}
|E(t,j) - E(t,0)| \le  \frac{d_{j+1}}{W}e^{-\int\limits_0^t
\gamma_{**}(\tau)\, d\tau}, \label{f13}
\end{equation}
\noindent gives us the corresponding speed of convergence to zero, for any $j$ and any $t \ge 0$.

\bigskip

Put now  $d_0=1$ and $d_{n+1}=(1+\varepsilon)d_n$ for $n \ge 0$, for a small positive $\varepsilon$.

\smallskip

Then also similar to \cite{Ammar2020}, we can obtain the following explicit bounds.

{\bf Proposition 1.}  Let there exist $\varepsilon>0$ such that
\begin{equation}
\int_0^\infty \left(\gamma^*(t)-\varepsilon \upsilon(t)\right)\, dt = +\infty, \label{f11*}
\end{equation}
\noindent where  $\upsilon(t) = \max\left(\eta(t), \lambda(t)\right)$.
Then:
\begin{eqnarray}
\left\|{\bf p}^{*}(t)-{\bf p}^{**}(t)\right\|_{1\textsf{D}} \le
e^{-\int\limits_0^t \left(\gamma^*(\tau)-\varepsilon \upsilon(\tau)\right)\, d\tau}\left\|{\bf
p}^{*}(0)-{\bf p}^{**}(0)\right\|_{1\textsf{D}}, \label{f12*}
\end{eqnarray}
\noindent and
\begin{equation}
|E(t,j) - E(t,0)| \le  \frac{d_{j+1}}{W}e^{-\int\limits_0^t
\left(\gamma^*(\tau)-\varepsilon \upsilon(\tau)\right)\, d\tau}, \label{f13*}
\end{equation}
\noindent for the corresponding $l_{1\textsf{D}}$ and $W$.

\bigskip

{\bf Second approach}, see also \cite{Zeifman2019amc}.

\bigskip Consider firstly the particular case of the same catastrophe rates, namely, suppose that all $\gamma_n(t)=\gamma^*(t)$. In this situation the equation (\ref{k1}) will look like this $r'(t) = -\eta(t) r(t) + \gamma^*(t)$, hence one can solve it:
\begin{equation}
r(t) = \int\limits_0^t e^{-\int\limits_\tau^t \eta(u)\, du}
\gamma^*(\tau) d\tau ,\label{r}
\end{equation}
\noindent because $r(0)=0$.

\smallskip

Consider  now the reduced forward Kolmogorov system (\ref{s2}) in
the form
\begin{equation}
\frac{d {\bf z} }{dt}=B\left( t\right) {{\bf z} }  +{\bf f}
\left(t\right), \quad t\ge 0, \label{s4}
\end{equation}
\noindent where  ${\bf f}\left(t\right)=\left(\eta(t)
r\left(t\right),0,0, \dots\right)^T$,  ${\bf
z}\left(t\right)=\left(p_0(t),p_1(t),\dots\right)^T$, and

$$B\left(t\right)=$$
{\tiny
    \begin{equation*}
    \left(
    \begin{array}{cccccccccccccc}
    \\
    -(\lambda (t)+\gamma^* (t))   & \mu(t) & 0 & 0 & \cdots & 0 & 0 & 0 & 0 & 0 & \cdots & 0  & 0 & \cdots \\
    \\
    \lambda (t) & -(\lambda (t)+\gamma^* (t)+\mu(t)) & \mu(t)  & 0 & \cdots & 0 & 0 & 0 & 0 & 0 & \cdots & 0  & 0 & \cdots\\

    \vdots & \vdots & \vdots & \vdots & \cdots & \vdots & \vdots & \vdots & \vdots & \vdots & \cdots & \vdots & \vdots & \cdots\\
    \\

    0 & 0 & 0 & 0 & \cdots & 0 & \lambda(t)  & -(\lambda (t) \beta (t)+\gamma^* (t)+\mu(t)) & \mu(t) &  0 & \cdots & 0 & 0 & \cdots\\
    \\
    0 & 0 & 0 & 0 & \cdots & 0 & 0 & \lambda(t) \beta (t)  & -(\lambda (t) \beta (t)+\gamma^* (t)+\mu(t)) & \mu(t) &  \cdots & 0 & 0 & \cdots\\
    \\
    \vdots & \vdots & \vdots & \vdots & \vdots & \vdots & \vdots & \vdots & \vdots & \vdots & \vdots & \vdots & \vdots & \vdots \\

    \end{array}
    \right).
    \label{m3}
    \end{equation*}
}

The solution of equation (\ref{s4}) can be written in the form
\begin{equation}
{\bf z}\left(t\right)= U_B\left(t,0\right){\bf
z}\left(0\right)+\int_0^t{}U_B\left(t,\tau\right){\bf
f}\left(\tau\right)\,d\tau, \label{f14}
\end{equation}
\noindent where  $U_B\left(t,s\right)$ is the Cauchy operator of the
corresponding homogeneous equation
\begin{equation}
\frac{d {\bf v} }{dt} = B\left( t\right) {{\bf v} }. \label{f15}
\end{equation}

\smallskip

Note that the  uniform estimate is completely analogous to Theorem
1, only with the replacement on the left side of ${\bf p}$ by ${\bf
z}$.

\bigskip

A significantly different situation with this approach arises when
we would like to consider general case, and even more to obtain weighted estimates.

\smallskip

Now we cannot find $r(t)$ in the closed form as in (\ref{r}). Instead of this put $r(t)=1-\sum_{i \ge 0} p_i(t)$. Then again we get
the equation (\ref{s4}) with another $B(t)$,
$$B\left(t\right)=$$
{\tiny
    \begin{equation*}
    \left(
    \begin{array}{cccccccccccccc}
    \\
    -(\lambda (t)+\gamma_0 (t)+\eta(t))   & \mu(t)-\eta(t) & -\eta(t) & -\eta(t) & \cdots & \cdots & \cdots & \cdots & \cdots &\cdots  & \cdots &\cdots   &\cdots  & \cdots \\
    \\
    \lambda (t) & -(\lambda (t)+\gamma_1 (t)+\mu(t)) & \mu(t)  & 0 & \cdots & 0 & 0 & 0 & 0 & 0 & \cdots & 0  & 0 & \cdots\\

    \vdots & \vdots & \vdots & \vdots & \cdots & \vdots & \vdots & \vdots & \vdots & \vdots & \cdots & \vdots & \vdots & \cdots\\
    \\

    0 & 0 & 0 & 0 & \cdots & 0 & \lambda(t)  & -(\lambda (t) \beta (t)+\gamma_k (t)+\mu(t)) & \mu(t) &  0 & \cdots & 0 & 0 & \cdots\\
    \\
    0 & 0 & 0 & 0 & \cdots & 0 & 0 & \lambda(t) \beta (t)  & -(\lambda (t) \beta (t)+\gamma_{k+1} (t)+\mu(t)) & \mu(t) &  \cdots & 0 & 0 & \cdots\\
    \\
    \vdots & \vdots & \vdots & \vdots & \vdots & \vdots & \vdots & \vdots & \vdots & \vdots & \vdots & \vdots & \vdots & \vdots \\

    \end{array}
    \right).
    \label{m300}
    \end{equation*}
}

Moreover, now  ${\bf f}\left(t\right)=\left(\eta(t),0,0, \dots\right)^T$, while ${\bf
z}\left(t\right)=\left(p_0(t),p_1(t),\dots\right)^T$.

$~~~~~~~$Let
\begin{eqnarray}
\mathfrak{D} =
\left(
\begin{array}{cccc}
d_0 & d_0 & d_0 & \cdots \\
\\
0 & d_1 & d_1 & \cdots  \\
\\
0 & 0 & d_2 & \cdots  \\
\\
& \ddots & \ddots & \ddots \\
\end{array}
\right)
\end{eqnarray}
and $ B^*(t) = \mathfrak{D}B(t)\mathfrak{D}^{-1}(t)= $
{\tiny
    \begin{equation}
    \left(
    \begin{array}{cccccccccccccc}
    \\
    -(\eta(t)+\gamma_0 (t))   & \frac{d_0}{d_1}(\gamma_0(t)-\gamma_1 (t)) & \frac{d_0}{d_2}(\gamma_1(t)-\gamma_2 (t)) & \cdots & \cdots & \cdots  \\
    \\
    \frac{d_1}{d_0} \lambda (t) & -(\lambda (t)+\gamma_1 (t)+\mu(t)) & \frac{d_1}{d_2}(\mu(t)+\gamma_1(t)-\gamma_2 (t))  & \frac{d_1}{d_3}(\gamma_2(t)-\gamma_3 (t)) & \cdots & \cdots\\

    \vdots & \vdots & \vdots & \vdots & \cdots & \vdots \\
    \\
    \vdots & \vdots & \vdots & \vdots & \vdots & \vdots  \\
    \end{array}
    \right).
    \label{m311}
    \end{equation}
}
Put
\begin{equation}
\gamma_B (t) =  \inf_i \left(|b^*_{ii}(t)| - \sum_{j\neq i} b^*_{ji}(t)\right).
\label{f18}
\end{equation}
\noindent  We have
\begin{eqnarray}
\|B^*(t)\| = \|B(t)\|_{1\mathfrak{D}} = \|\mathfrak{\mathfrak{D}} B(t) \mathfrak{D}^{-1}\| =
\sup_i \left(|b^*_{ii}(t)| +
\sum_{j\neq i} b^*_{ji}(t)\right) = \qquad \qquad \\
= \sup_i \left(2|b^*_{ii}(t)| + \sum_{j\neq i}
b^*_{ji}(t)- |b^*_{ii}(t)|\right) \le  2 \sup_i
|b^*_{ii}(t)| - \gamma_B(t)\le \ 2L -\gamma_B(t), \nonumber
\label{bound302}
\end{eqnarray}
\noindent  hence the operator function $B(t)$ is bounded on the space $l_{1\mathfrak{D}}$. Therefore we can apply the same approach to equation (\ref{s4}) in the space $l_{1\mathfrak{D}}$. Now the equality
\begin{eqnarray} \gamma
\left(B(t)\right)_{1\mathfrak{D}} = \gamma \left(\mathfrak{D} B(t)
\mathfrak{D}^{-1}\right) =  \sup_i \left(b^*_{ii}(t) + \sum_{j\neq
    i}b^*_{ji}(t)\right) = - \gamma_B(t),
\label{f19}
\end{eqnarray}
implies the following statement.

\smallskip

{\bf Theorem 3.} \label{t4} Let  
\begin{equation}
\int_0^\infty \gamma_B(t)\, dt = +\infty,
\label{f20}
\end{equation}
\noindent for some $\{d_i\}$, $1 = d_0 \le d_1 \le
\dots$.
Then  $X(t)$ is weakly ergodic and
\begin{eqnarray}
\left\|{z}^{*}(t)-{z}^{**}(t)\right\|_{1\textsf{D}} \le e^{-\int\limits_0^t \gamma_B(\tau)\,
    d\tau}\left\|{z}^{*}(0)-{z}^{**}(0)\right\|_{1\textsf{D}},
\label{f21}
\end{eqnarray}
\noindent for any initial conditions ${z}^{*}(0), {z}^{**}(0)$ and any $t \ge 0$.

\smallskip

Put now $W = \inf_{k \ge 1}\frac{d_k}{k}$. Then $W\|{\bf z}\|_{1E} \le \|{\bf z}\|_{1\textsf{D}}$.

\medskip

{\bf Corollary 2.}   \label{cor2} Let under assumptions of Theorem 3, in addition $W > 0$. Then   the following bound holds
\begin{equation}
|E(t,j) - E(t,0)| \le  \frac{1+d_j}{W}e^{-\int\limits_0^t
    \gamma_B(\tau)\, d\tau}, \label{3012e}
\end{equation}
for any $j$ and any $t \ge 0$.

\bigskip

{\bf Remark.} One can put $d_0=1$, $d_1=\epsilon$, $d_{k+1}=(1+\epsilon)d_k$ for $k \ge 1$,
and obtain the analogue of Proposition 1 for the second approach.

\bigskip

{\bf Remark.} In all our statements, we can replace the condition of
monotonicity of the sequence $\{d_k\}$  by condition $d = \inf_k d_k
>0$, with the corresponding change in the estimates; see, for
example, \cite{Zeifman2020amcs}.

%=============================================================================================
\section{\textit{Perturbation bounds}}

Consider here the application of general perturbation bounding (see
the recent review in \cite{Zeifman2020mdpi}) for the models under
study.
\smallskip
Consider a   "perturbed" queue-length process$\bar X(t), t \ge 0$
with the corresponding transposed intensity matrix $\bar A(t)$,
where the "perturbation" matrix $\hat A(t) = A(t) -\bar A(t)$ is
small in a sense. Namely, we assume that the perturbed queue is of
the same nature as the original one. Hence, the perturbed intensity
matrix also has the same structure, with the corresponding perturbed
intensities  $\bar \eta(t)$,  $\bar \gamma_n(t)$, $\bar \lambda(t)$,
$\bar \mu(t)$,  $\bar \beta(t)$.
\smallskip
Let
\begin{eqnarray} \label{stab01} |\eta(t) - \bar\eta(t)| = |\hat\eta(t)| \le \hat\epsilon, \quad
|\gamma_n(t) - \bar\gamma_n(t)| =|\hat\gamma_n(t)|\le \hat\epsilon, \\
|\lambda(t) - \bar\lambda(t)| =|\hat\lambda(t)| \le \hat\epsilon, \quad
|\mu(t) - \bar\mu(t)| = |\hat\mu(t)| \le \hat\epsilon, \quad |\beta(t)
- \bar\beta(t)| = |\hat\beta(t)| \le \hat\epsilon. \nonumber
\end{eqnarray}
Hence
\begin{eqnarray} \label{stab01.5} |\lambda(t)\beta(t) - \bar\lambda(t) \bar\beta(t)|\le
\lambda(t)|\hat\beta(t)|+\bar\beta(t)|\hat\lambda(t)| \le (L+1)
\hat\epsilon.
\end{eqnarray}
 Then we obtain from (\ref{m1}) the following bound
\begin{eqnarray} \label{stab02} \|\hat A(t) \| = 2\sup_{k}\left|\hat
a_{kk}\left(t\right)\right| =\\
2\max\left(|\hat\eta(t)|,|\hat\lambda(t)|+
|\hat\gamma(t)|,|\hat\lambda(t)|+ |\hat\gamma_n(t)|+|\hat\mu(t)|,(L+1)
\hat\epsilon + |\hat\gamma(t)|+|\hat\mu(t)|\right) \le
\left(2L+6\right)\hat\epsilon.\nonumber
\end{eqnarray}

\medskip

Firstly we formulate the  perturbation bounds for the vector of state probability in the situation of Theorem 1.

The next statement follows immediately from Theorem 1
\cite{Zeifman2020mdpi} (see also the first corresponding homogeneous
result in \cite{Mitrophanov2003} and for inhomogeneous situation in
\cite{Zeifman2012}).

{\bf Theorem 4.} \label{t4} Let under assumption of Theorem 1 the
catastrophe intensity $\gamma(t)$ be such that
\begin{equation}
e^{-\int\limits_s^t \gamma^*\left(\tau\right)\, d\tau} \le N e^{-\gamma_0(t-s)},
 \label{stab03}\end{equation}
 \noindent for some positive $N, \gamma_0$. Then  the following perturbation bound holds:
\begin{eqnarray}
\limsup_{t\to \infty}  \|{\bf p}(t) - \bar {\bf p}(t)\| \le
\frac{\hat\epsilon\left(2L+6\right)\left(
1+\log{(N/2)}\right)}{\gamma_0},
 \label{stab04} \end{eqnarray}
\noindent for any perturbed queue with the respectively closed intensities satisfying to (\ref{stab01}).

\bigskip

Stability bound from Theorem 2 is based on results of \cite{Zeifman2015i,Zeifman2016ipiran}.

Note that (\ref{m2}), (\ref{stab01}) and (\ref{stab02}) imply the
inequality:
\begin{eqnarray} \label{stab11} \|\hat A^*(t) \| = 2\sup_{k}\left|\hat
a_{kk}^*\left(t\right)\right|  \le \|\hat A(t) \| \le
\left(2L+6\right)\hat\epsilon.
\end{eqnarray}
\smallskip

On the other hand, we have $\|{\bf g}(t)\|_{1D} = \gamma^*(t) \le L$
for almost all $t \ge 0$.

\smallskip

Then  Theorem 4 from \cite{Zeifman2016ipiran} imply the next
statement.

\bigskip

{\bf Theorem 5.} \label{t6} Let under assumptions of Theorem 2 the
following estimates hold:
\begin{equation}
e^{-\int\limits_s^t \gamma_{**}\left(\tau\right)\, d\tau} \le N^{**}
e^{-\gamma_0^{**}(t-s)},
 \label{stab12}\end{equation}
 \noindent for some positive $N^{**}, \gamma_0^{**}$, and
\begin{equation}
H = \sup_{|i-j|=1}\frac{d_i}{d_j} < \infty.
 \label{stab13}\end{equation}
Then
\begin{eqnarray}
\limsup_{t\to \infty}  \|{\bf p}(t) - \bar {\bf p}(t)\|_{1D} \le
\frac{\left(4L+12\right) \hat\epsilon HL
(N^{**})^2}{\gamma_0^{**}\left(\gamma_0^{**} -
\left(4L+12\right)\hat\epsilon H\right)}. \label{stab14}
\end{eqnarray}
Moreover, if $W = \inf_{k \ge 1}\frac{d_{k+1}}{k} >0$, then
\begin{eqnarray}
\limsup_{t\to \infty}  |E(t,0) - \bar E(t,0)| \le
\frac{\left(4L+12\right) \hat\epsilon HL
(N^{**})^2}{W\gamma_0^{**}\left(\gamma_0^{**} -
\left(4L+12\right)\hat\epsilon H\right)}. \label{stab15}
\end{eqnarray}

\bigskip

\bigskip

Finally, we obtain perturbation bounds based on the ergodicity
estimates of Theorem 3.

\bigskip

{\bf Theorem 6.} \label{t5} Let under assumptions of Theorem 3
the
following estimates hold:
\begin{equation}
e^{-\int\limits_s^t \gamma_{B}\left(\tau\right)\, d\tau} \le N^{B}
e^{-\gamma_0^{B}(t-s)},
 \label{stab22}\end{equation}
 \noindent for some positive $N^{B}, \gamma_0^{B}$, and
\begin{equation}
H = \sup_{|i-j|=1}\frac{d_i}{d_j} < \infty.
 \label{stab23}\end{equation}
Then
\begin{eqnarray}
\limsup_{t\to \infty}  \|{\bf p}(t) - \bar {\bf
p}(t)\|_{1\textsf{D}} \le \frac{\hat\epsilon N^{B}
\left(L+1\right)\left(6HL
N^{B}+\gamma_0^{B}\right)}{\gamma_0^{B}\left(\gamma_0^{B} -
12\hat\epsilon H N^{B}\left(L+1\right)\right)}. \label{stab24}
\end{eqnarray}
Moreover, if $W = \inf_{k \ge 1}\frac{d_{k}}{k} >0$, then
\begin{eqnarray}
\limsup_{t\to \infty}  |E(t,0) - \bar E(t,0)| \le \frac{\hat\epsilon
N^{B} \left(L+1\right)\left(6HL
N^{B}+\gamma_0^{B}\right)}{W\gamma_0^{B}\left(\gamma_0^{B} -
12\hat\epsilon H N^{B}\left(L+1\right)\right)}. \label{stab25}
\end{eqnarray}

{\bf Proof.} It is sufficient to note that
\begin{equation} \|B(t) - \bar
B(t)\|_{1\textsf{D}} \le H \|B(t) - \bar B(t)\|_{1}  \le H \|A(t) -
\bar A(t)\|_{1}  \le \left(2L+6\right)H \hat\epsilon
 \label{stab20}\end{equation}
\noindent and
\begin{equation} \|{\bf f}(t) - \bar {\bf f}(t)\|_{1\textsf{D}} = ||\eta(t) r(t) - \bar\eta(t)\bar r(t)| \le
\left(L+1\right) \hat\epsilon.
 \label{stab21}\end{equation}

Then our claim follows from  Theorem 2 \cite{Zeifman2020mdpi}.

%=============================================================================================
\section{\textit{Numerical Examples}}

In this section, we will review  two numerical examples to support the results obtained as well as clarify the nature of the behavior of the proposed system. Where through these examples, we will apply the analytical results obtained in the theories corollaries in the previous sections. For these examples, we would consider the threshold value equal to 100 ($k=100$).

{\bf Example 1.} Let our queueing model have  the following rates:
\noindent $\eta(t)= 3 +\sin 2 \pi t$, $\gamma_n(t)=2+0.5\cos 2 \pi t$, for any $n$,
$\lambda(t)= 10+10\sin 2 \pi t$, $\mu(t)=2+\cos2 \pi t$,
$\beta(t)=0.7$.

\smallskip

Apply all our bounds for this specific situation.

For Theorem 1 and the respective "stability" Theorem 4  we need $L$, $N$ and $\gamma_0$. Obviously we have $L \le 25.5$, and $\gamma^*(t)= 2+0.5\cos 2 \pi t$. Consider now
$$e^{-\int\limits_s^t \gamma^*\left(\tau\right)\, d\tau} = e^{-2(t-s) -\frac{\sin 2 \pi t - \sin 2 \pi s}{4 \pi}}  \le e^{-2(t-s) +\frac{1}{2 \pi}} \le 2 e^{-2(t-s)},$$
hence one can put $N=2$ and $\gamma_0=2$ in  (\ref{stab03}).

\smallskip

For applying of Theorems 2, 5 we put $\varepsilon = 0.05$, $d_0=1$
and $d_{k+1}=(1+\varepsilon)d_k$ for $k \ge 0$. Then have
$H=1+\varepsilon<2$, $\upsilon(t) = \max\left(\eta(t),
\lambda(t)\right) = 10+10\sin 2 \pi t$, and $\gamma_{**}(t)=
\gamma_{\varepsilon}(t) =1.5+ 0.5\cos 2 \pi t-0.5\sin 2 \pi t$.
$$ e^{-\int\limits_s^t \gamma_{**}\left(\tau\right)\,
d\tau}   \le e^{-1.5(t-s) +\frac{1}{4 \pi}} \le 2 e^{-1.5(t-s)},$$
therefore one can get $N^{**}=2$ and $ \gamma_0^{**}=1.5$ in
(\ref{stab12}).

\smallskip

Finally, for applying of Theorems 2, 5 we put $\varepsilon = 0.05$,
$d_0=1$,  $d_1=\varepsilon$, and $d_{k+1}=(1+\varepsilon)d_k$ for $k
\ge 1$. Then we have $H=\frac{1}{\varepsilon}$, and
$\gamma_{B}(t)=\gamma_{**}(t)= \gamma_{\varepsilon}(t)$, hence one
can put $N^{B}=2$, $\gamma_0^{B}=1.5$ in (\ref{stab22}).

\bigskip

Now we have from Theorem 1:

\smallskip

\begin{eqnarray}
\left\|{{\bf p}}^{*}\left(t\right)-{{\bf
p}}^{**}\left(t\right)\right\|\le 4 e^{-2 t }\left\|{{\bf
p}}^{*}\left(0\right)-{{\bf p}}^{**}\left(0\right)\right\|;
\label{th1_1}
\end{eqnarray}
\smallskip
\noindent and 
\begin{eqnarray}
\left\|{\bf p}^{*}(t)-{\bf p}^{**}(t)\right\|_{1\textsf{D}} \le 2
e^{-1.5t }\left\|{\bf p}^{*}(0)-{\bf
p}^{**}(0)\right\|_{1\textsf{D}}, \label{th2_1a}
\end{eqnarray}
\begin{equation}
|E(t,j) - E(t,0)| \le  \frac{2\cdot{1.05}^{j+1}}{W}e^{-1.5 t },
\label{th2_1b}
\end{equation}
\noindent from Theorem 2 and Corollary 1, and almost the same from
Theorem 3 and Corollary 2.

\medskip

The corresponding perturbation bounds are:

\smallskip

\begin{eqnarray}
\limsup_{t\to \infty}  \|{\bf p}(t) - \bar {\bf p}(t)\| \le 30
\hat\epsilon,
 \label{th4_1} \end{eqnarray}
\noindent from Theorem 4;

\smallskip

\begin{eqnarray}
\limsup_{t\to \infty}  \|{\bf p}(t) - \bar {\bf p}(t)\|_{1D} \le
\frac{4\cdot 10^5 \hat\epsilon }{1 - 4\cdot 10^3\hat\epsilon }, \label{th5_1a}
\end{eqnarray}
\noindent and
\begin{eqnarray}
\limsup_{t\to \infty}  |E(t,0) - \bar E(t,0)| \le \frac{4 \cdot 10^5
\hat\epsilon }{W\left(1 - 4 \cdot 10^3\hat\epsilon\right)}, \label{th5_1b}
\end{eqnarray}
\noindent from Theorem 5; and bounds of Theorem 6 are much worse.

\bigskip

One can note that for this model bounds from Theorems 3 and 6 are
worse because of the matrices $B(t)$ and $B^*(t)$ have very special
structure.

\bigskip

Further we apply the approximations by truncated processes, see details in 
\cite{Zeifman2014i,Zeifman2017tpa}. Namely, it is sufficient to put the
dimensionality of the truncated process $200$ and  the
corresponding time interval $[0,20]$, hence the limit interval itself is $[19,20]$. 
Now Figures 1--4 shows us the behavior of the probability of the empty queue and the mean
respectively. In Figures 5-6 one can see the  perturbation bounds
for the corresponding limiting characteristics with $\hat\epsilon
=10^{-3}$  for bound (\ref{th4_1}) and $\hat\epsilon =10^{-6}$
 for (\ref{th5_1b}).

\bigskip

\begin{figure}
\centering
\includegraphics[height=20cm]{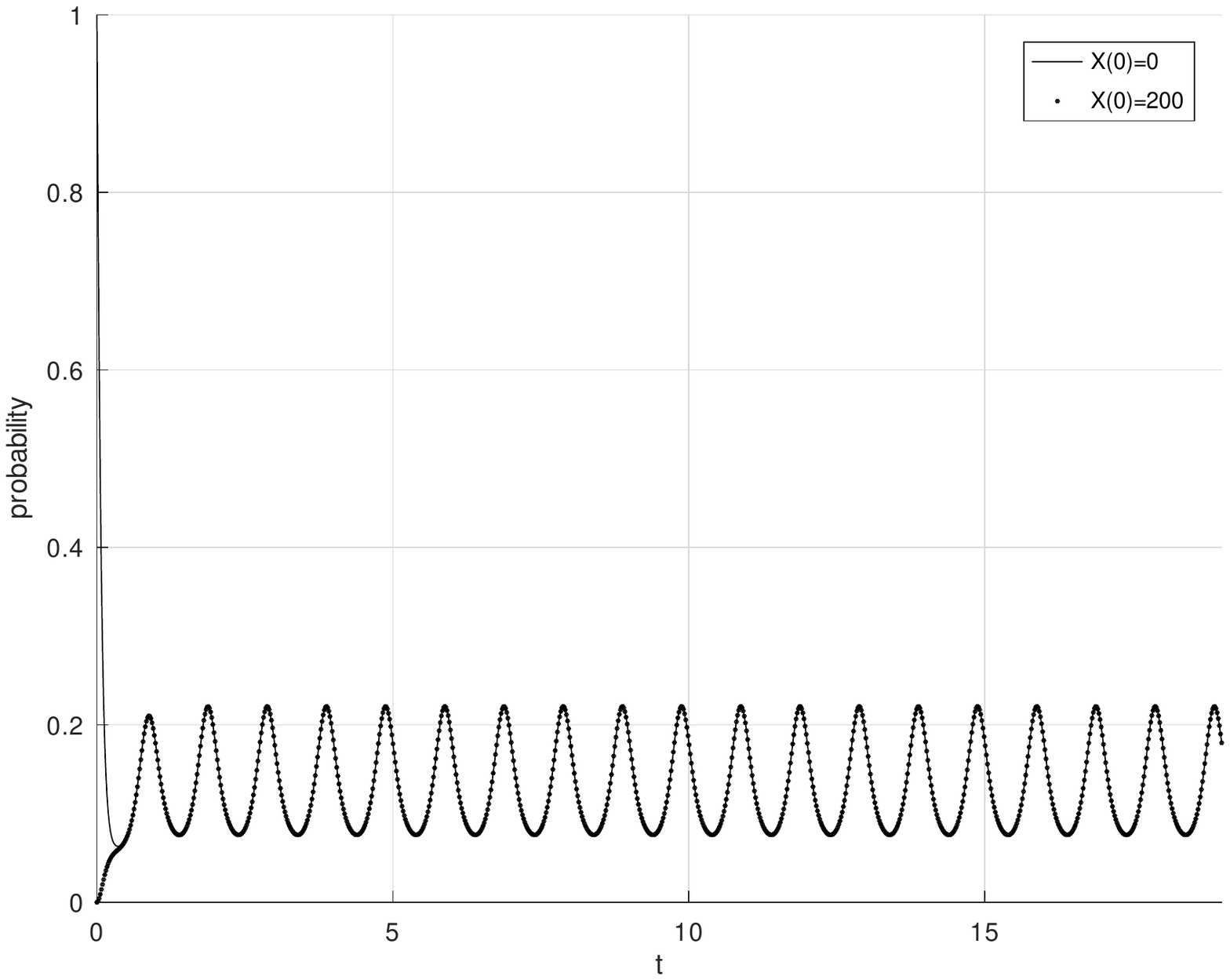}
\vspace{-5cm}\caption{Example 1. Probability of the empty queue for $ t \in
[0,19]$.}
\end{figure}

\begin{figure}
\centering
\includegraphics[height=20cm]{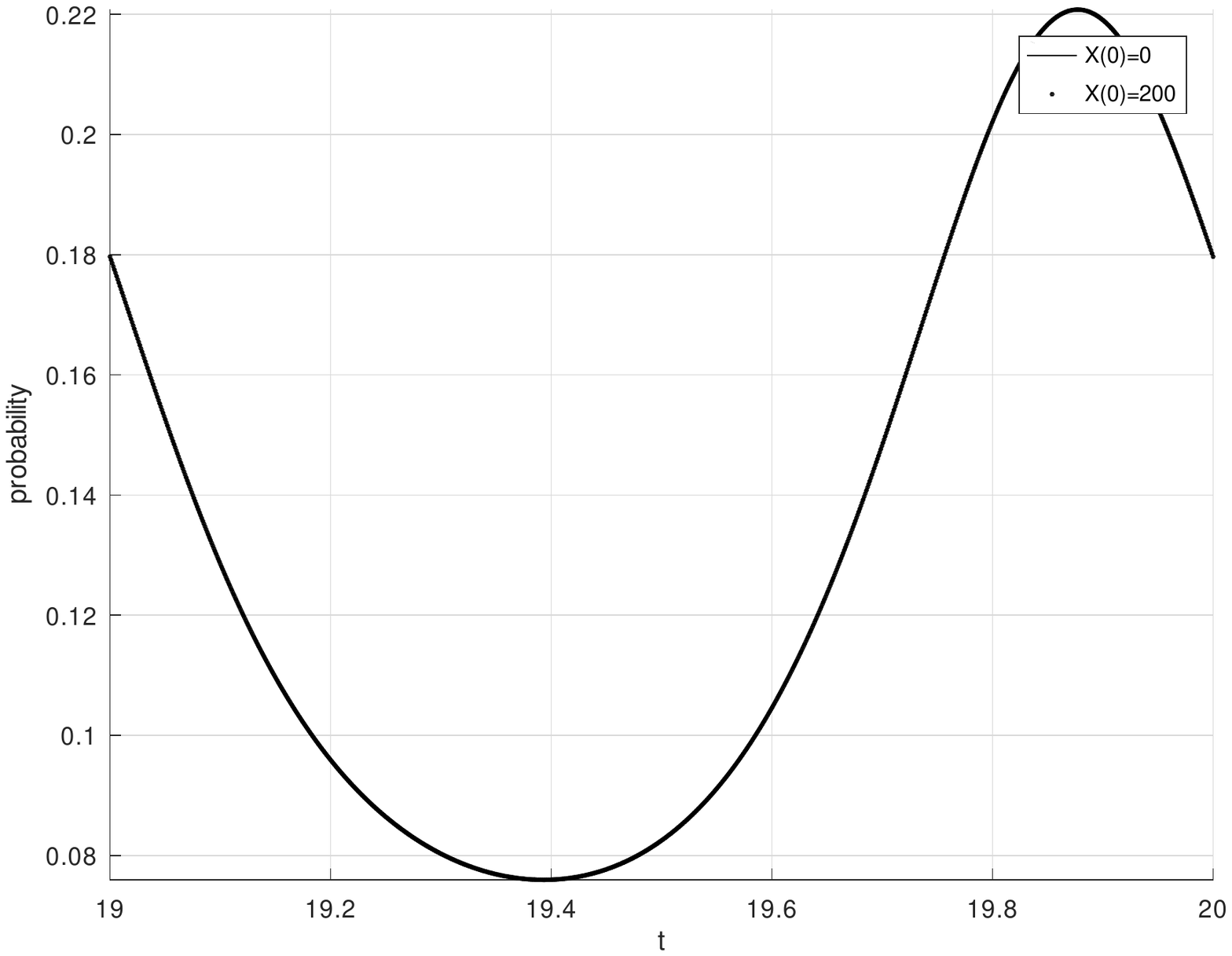}
\vspace{-5cm}\caption{Example 1. Approximation of the limiting probability of
 empty queue for $ t \in [19,20]$.}
\end{figure}

\begin{figure}
\centering
\includegraphics[height=20cm]{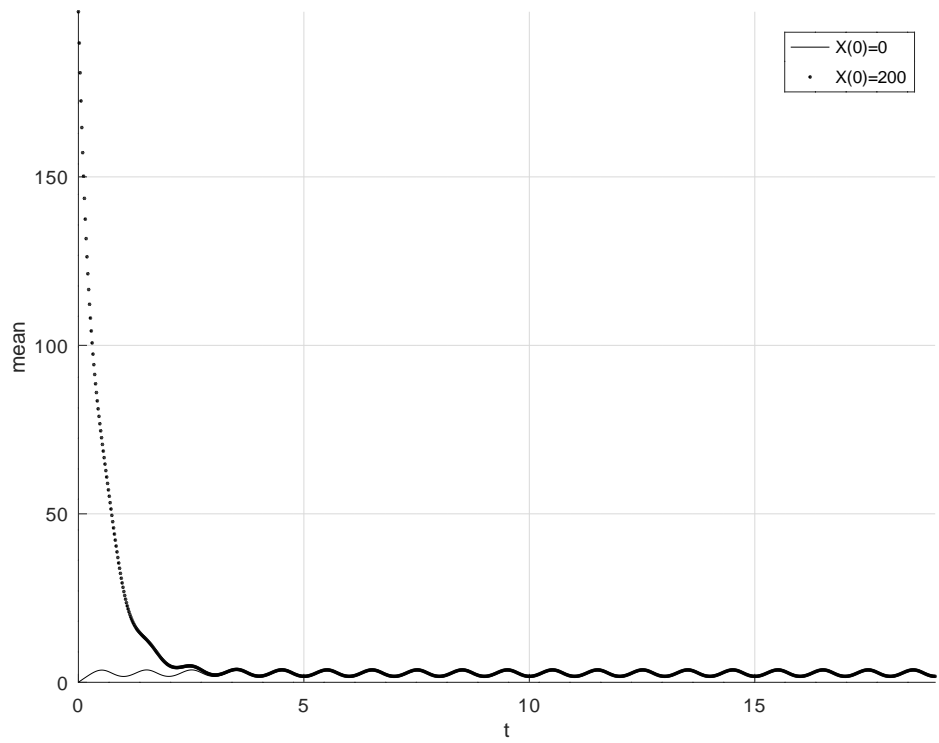}
\vspace{-5cm} \caption{Example 1. The mean $E(t,k)$ for $ t \in [0,19]$.}
\end{figure}

\begin{figure}
\centering
\includegraphics[height=20cm]{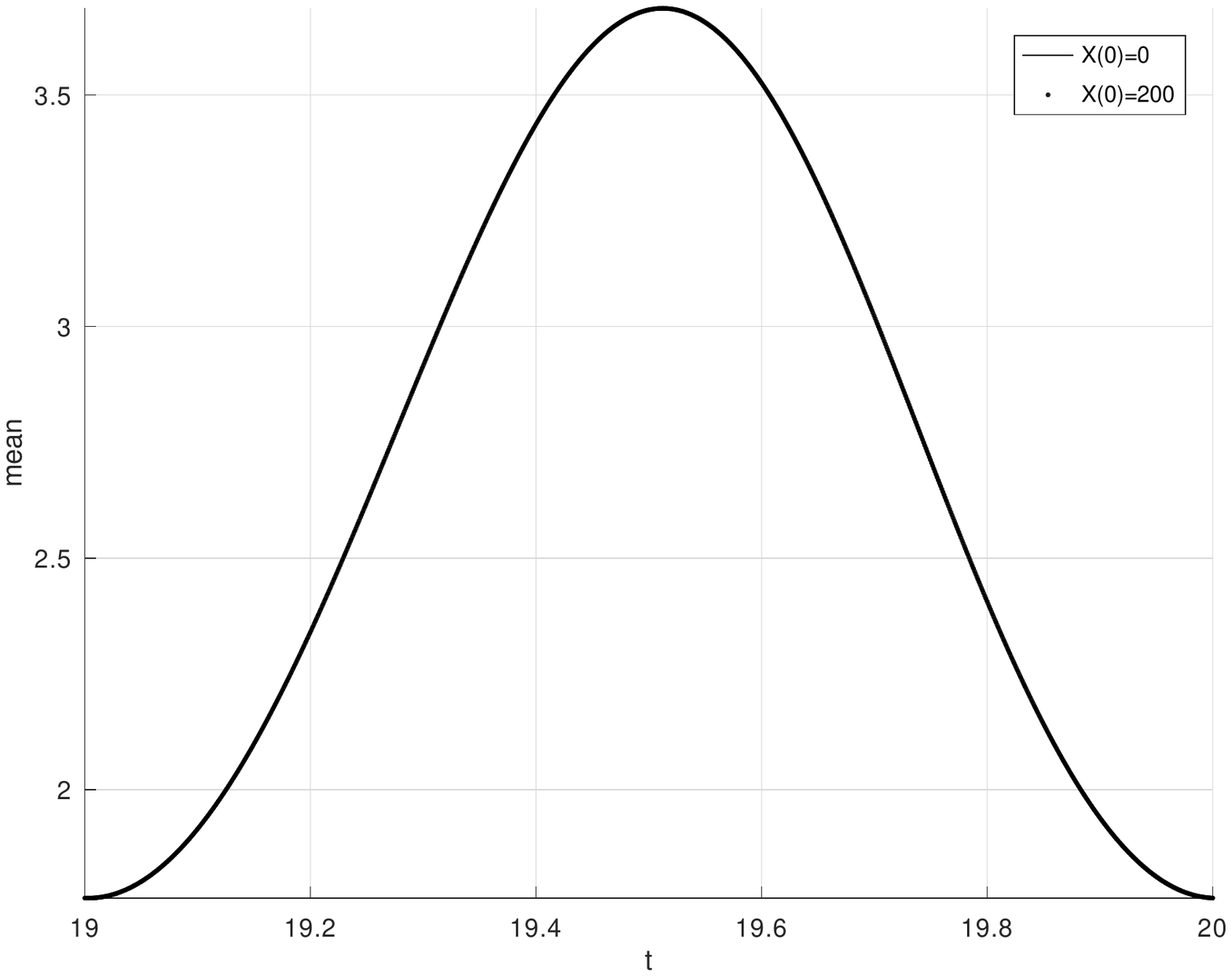}
\vspace{-3cm}\caption{Example 1. Approximation of the limiting mean $E(t,k)$
for $ t \in [19,20]$.}
\end{figure}

\begin{figure}
\centering
\includegraphics[height=20cm]{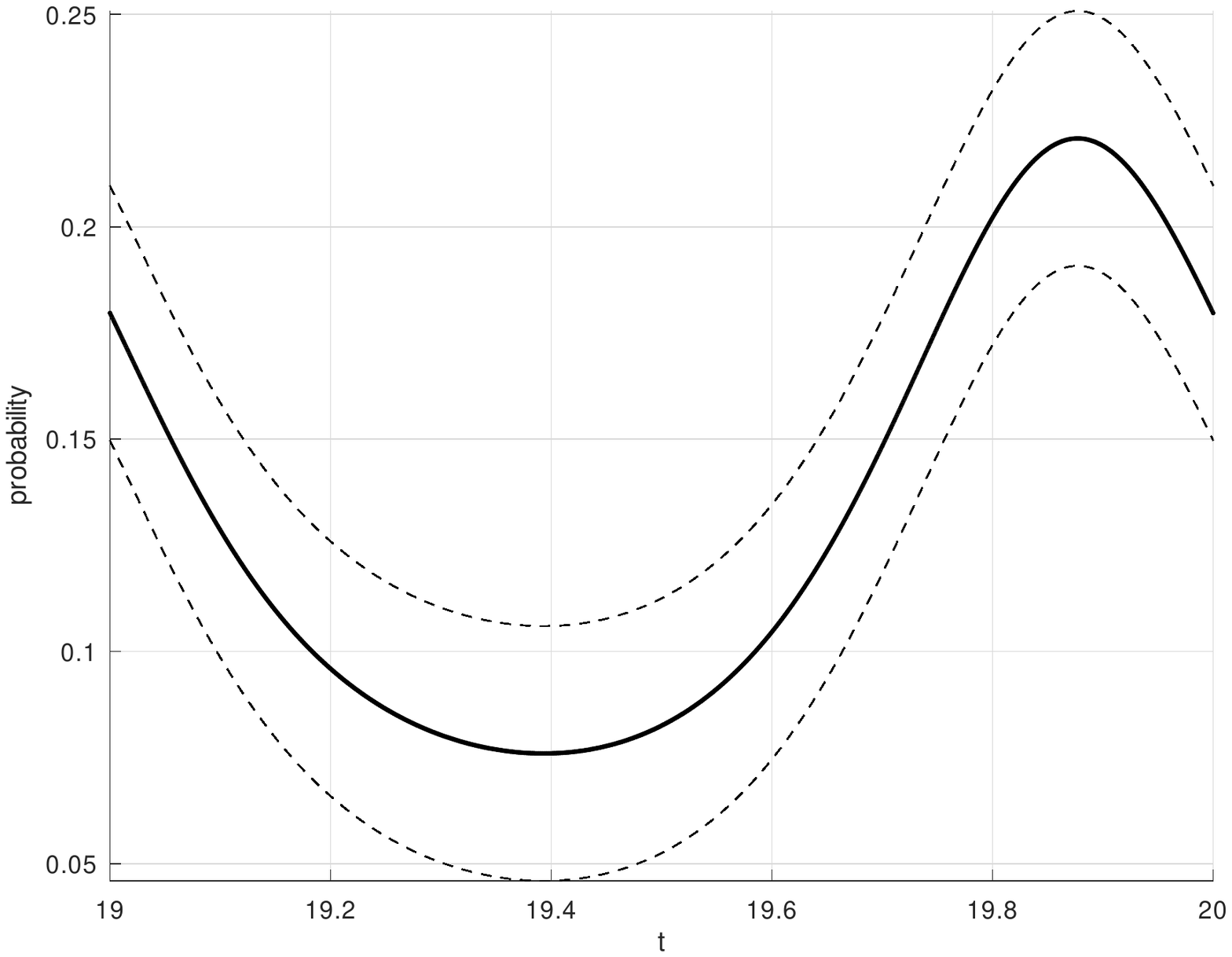}
\vspace{-5cm}\caption{Example 1. Perturbation bounds for the limiting
probability of  empty queue for $ t \in [19,20]$.}
\end{figure}

\begin{figure}
\centering
\includegraphics[height=20cm]{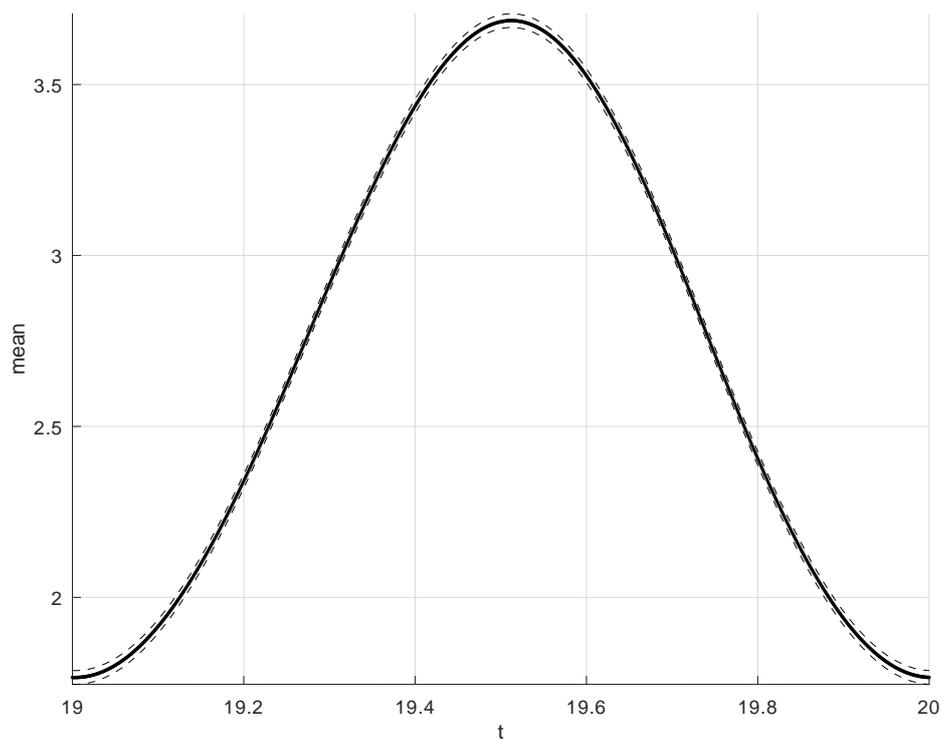}
\vspace{-5cm}\caption{Example 1. Perturbation bounds for  the limiting mean
$E(t,0)$ for $ t \in [19,20]$.}
\end{figure}

\bigskip

{\bf Example 2.} Consider now the  model with the following rates:
\noindent $\eta(t)= 3 +\sin 2 \pi t$, $\gamma_0(t)=2+0.5\cos 2 \pi t$, $\gamma_k(t)=0$ if $k \ge 1$,
$\lambda(t)= 1+\sin 2 \pi t$, $\mu(t)=5+\cos2 \pi t$, $\beta(t)=0.7$.

\smallskip

Here we have $\gamma^* =0$ and the first approach is not applicable.

For using Theorems 3 and 6 we firstly put $d_0=1$, $d_1=2.5$, and $d_{k+1}=(1.5)d_k$ for $k \ge 1$.

Then we have:
$$ |b^*_{00}(t)| - \sum_{j\neq 0} b^*_{j0}(t) = \eta(t)+\gamma_0(t) -  \frac{5}{2}\lambda (t) \ge   0.5,$$

$$ |b^*_{11}(t)| - \sum_{j\neq 1} b^*_{j1}(t) =  \mu(t) - 0.5\lambda (t)-\frac{2}{5}\gamma_0(t) \ge 2,$$

$$ |b^*_{ii}(t)| - \sum_{j\neq i} b^*_{ji}(t) = \frac{0.5}{1.5}\mu(t) - 0.5\lambda (t) \ge \frac{1}{3},
\quad 2 \le i \le k - 2,$$

$$ |b^*_{ii}(t)| - \sum_{j\neq i} b^*_{ji}(t) =
\frac{0.5}{1.5}\mu(t) + \lambda(t) (1 - \beta(t)) - 0.5\lambda (t) \beta(t) \ge \frac{1}{3},
\quad i =k-1,$$

$$ |b^*_{ii}(t)| - \sum_{j\neq i} b^*_{ji}(t) =
\frac{0.5}{1.5}\mu(t) - 0.5\lambda (t) \beta(t) \ge \frac{1}{3},
\quad i \ge k.$$

Hence,  in (\ref{f19}) we get
$$ \gamma_{B}(t) =  \inf_i \left(|b^*_{ii}(t)| - \sum_{j\neq i}
b^*_{ji}(t)\right) \ge \frac{1}{3},$$
\noindent and we obtain instead of (\ref{f21}) and (\ref{3012e}) the following bounds:
\begin{eqnarray}
\left\|{z}^{*}(t)-{z}^{**}(t)\right\|_{1\textsf{D}} \le e^{-\frac{t}{3}}\left\|{z}^{*}(0)-{z}^{**}(0)\right\|_{1\textsf{D}},
\label{ex2_01}
\end{eqnarray}
\noindent and
\begin{equation}
|E(t,j) - E(t,0)| \le  \frac{1+d_j}{W}e^{-\frac{t}{3}}. \label{ex2_02}
\end{equation}
Moreover, Theorem 6 gives us the corresponding perturbation bounds, namely, we have
$L= 8$, $H=2$, $N^{B}=1$, $\gamma_0^{B}=\frac{1}{3}$, and instead of (\ref{stab24}) and (\ref{stab25}) the following inequalities hold:
Then
\begin{eqnarray}
\limsup_{t\to \infty}  \|{\bf p}(t) - \bar {\bf
p}(t)\|_{1\textsf{D}} \le 10^4 \cdot \hat\epsilon, \label{ex2_03}
\end{eqnarray}
\noindent and
\begin{eqnarray}
\limsup_{t\to \infty}  |E(t,0) - \bar E(t,0)| \le 2\cdot 10^4 \cdot \hat\epsilon, \label{ex2_04}
\end{eqnarray}
\noindent for sufficiently small $\hat\epsilon >0$.

Further we compute the corresponding limiting characteristics. Figures 7--10 shows us
the behavior of the probability of the empty queue and the mean
respectively, and in Figures 11-12 one can see the  perturbation bounds
for the corresponding limiting characteristics with $\hat\epsilon
=10^{-4}$  for bound (\ref{ex2_03}) and $\hat\epsilon =5\cdot 10^{-5}$
 for (\ref{ex2_04}).

\bigskip

\begin{figure}
\centering
\includegraphics[height=20cm]{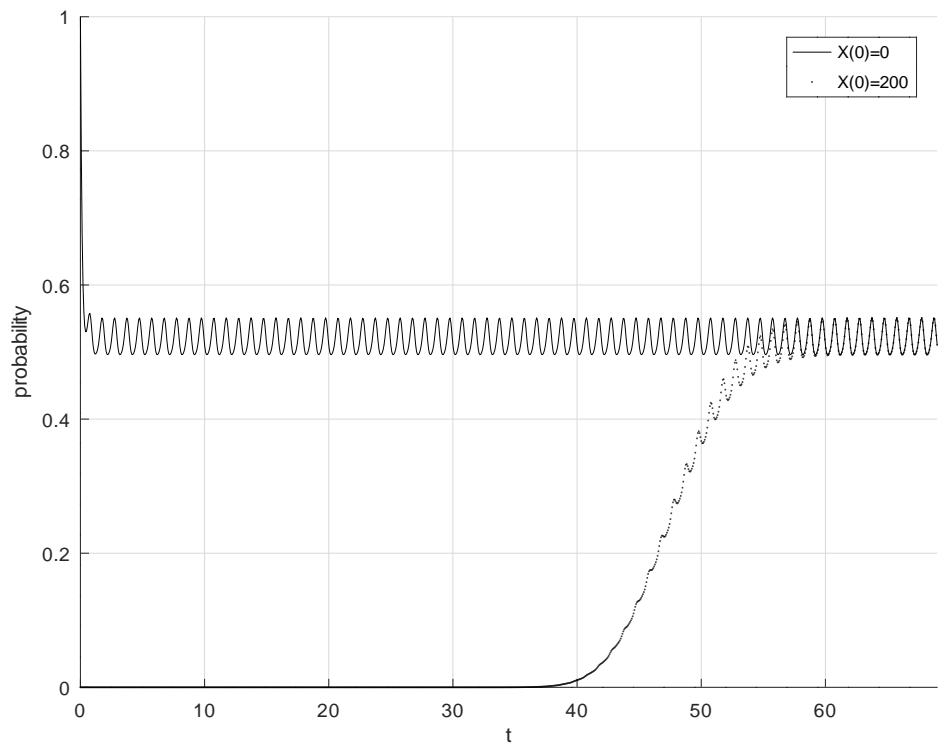}
\vspace{-5cm}\caption{Example 2. Probability of the empty queue for $ t \in
[0,69]$.}
\end{figure}

\begin{figure}
\centering
\includegraphics[height=20cm]{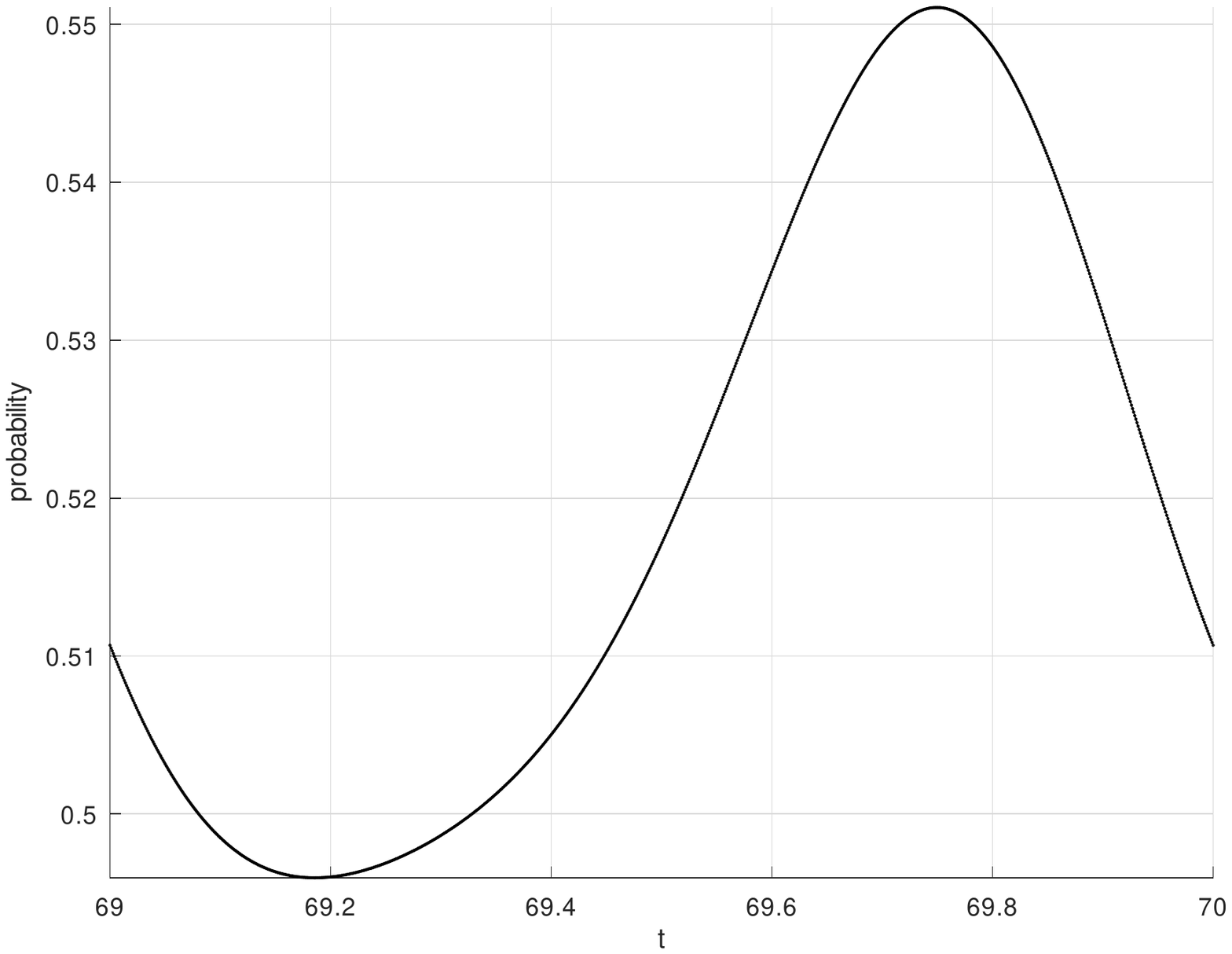}
\vspace{-5cm}\caption{Example 2. Approximation of the limiting probability of
 empty queue for $ t \in [69,70]$.}
\end{figure}

\begin{figure}
\centering
\includegraphics[height=20cm]{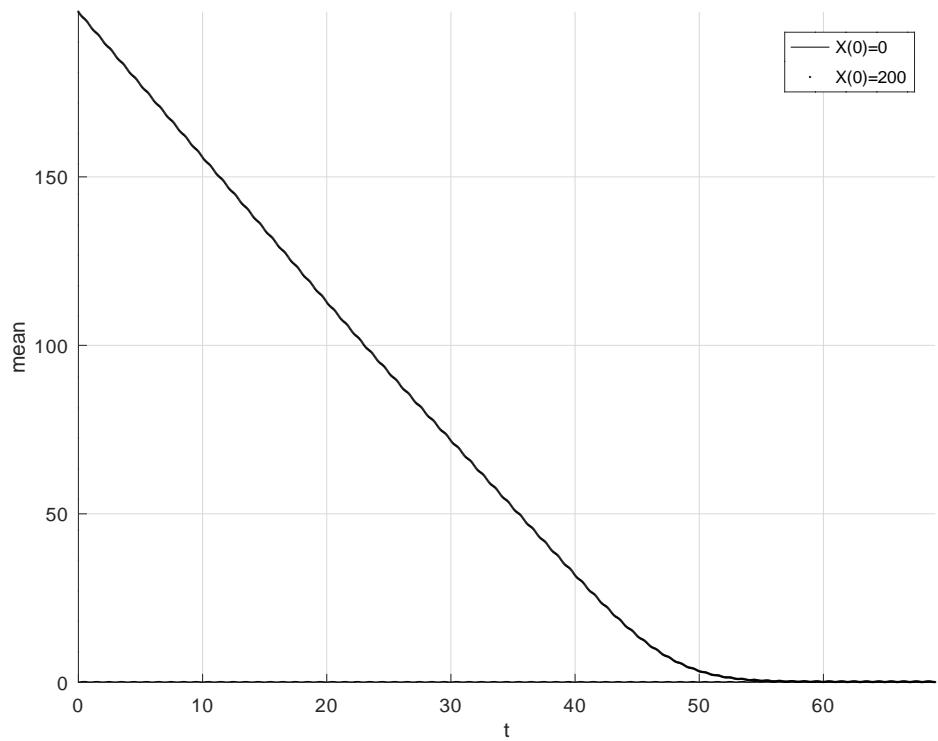}
\vspace{-5cm} \caption{Example 2. The mean $E(t,k)$ for $ t \in [0,69]$.}
\end{figure}

\begin{figure}
\centering
\includegraphics[height=20cm]{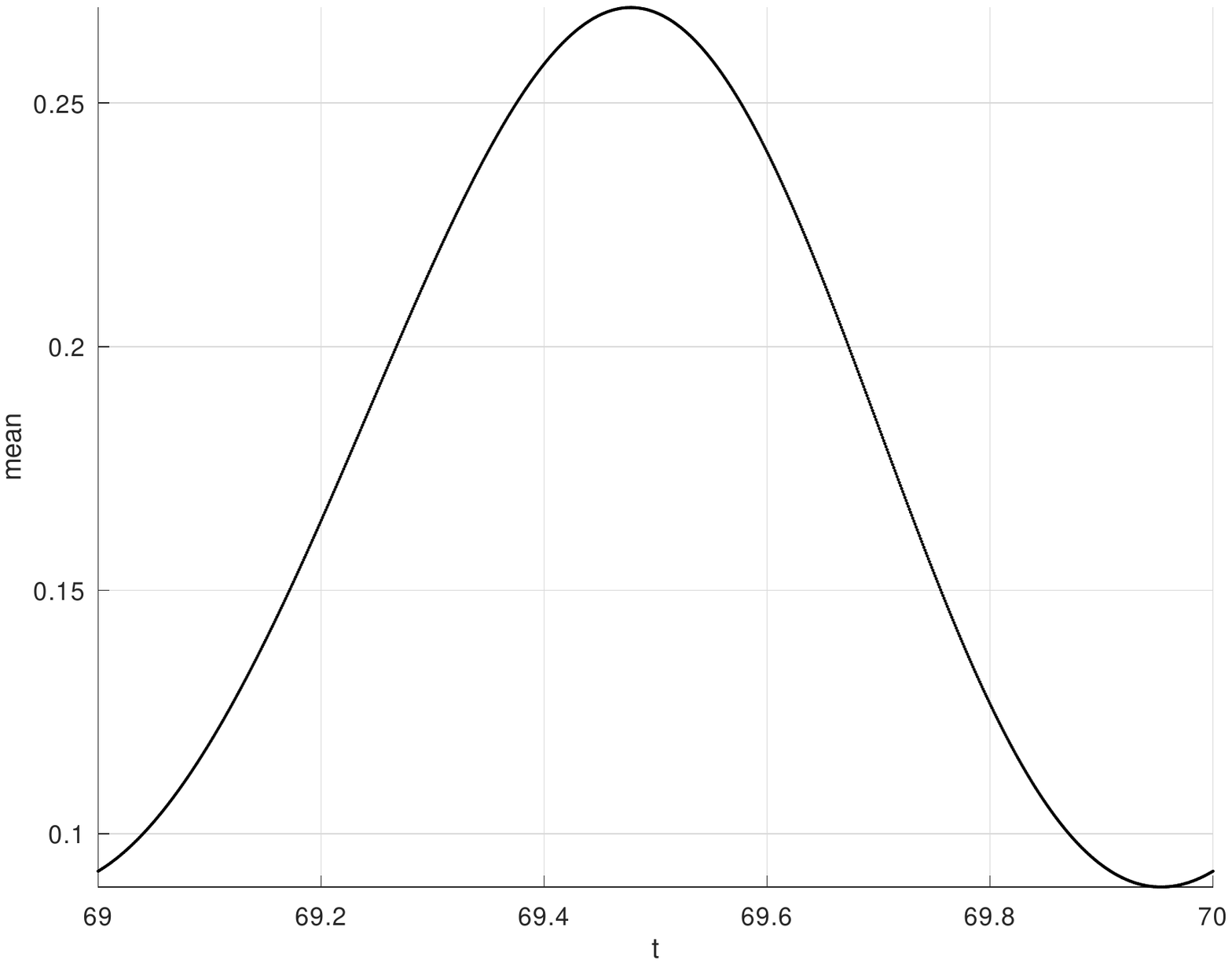}
\vspace{-3cm}\caption{Example 2. Approximation of the limiting mean $E(t,k)$
for $ t \in [69,70]$.}
\end{figure}

\begin{figure}
\centering
\includegraphics[height=20cm]{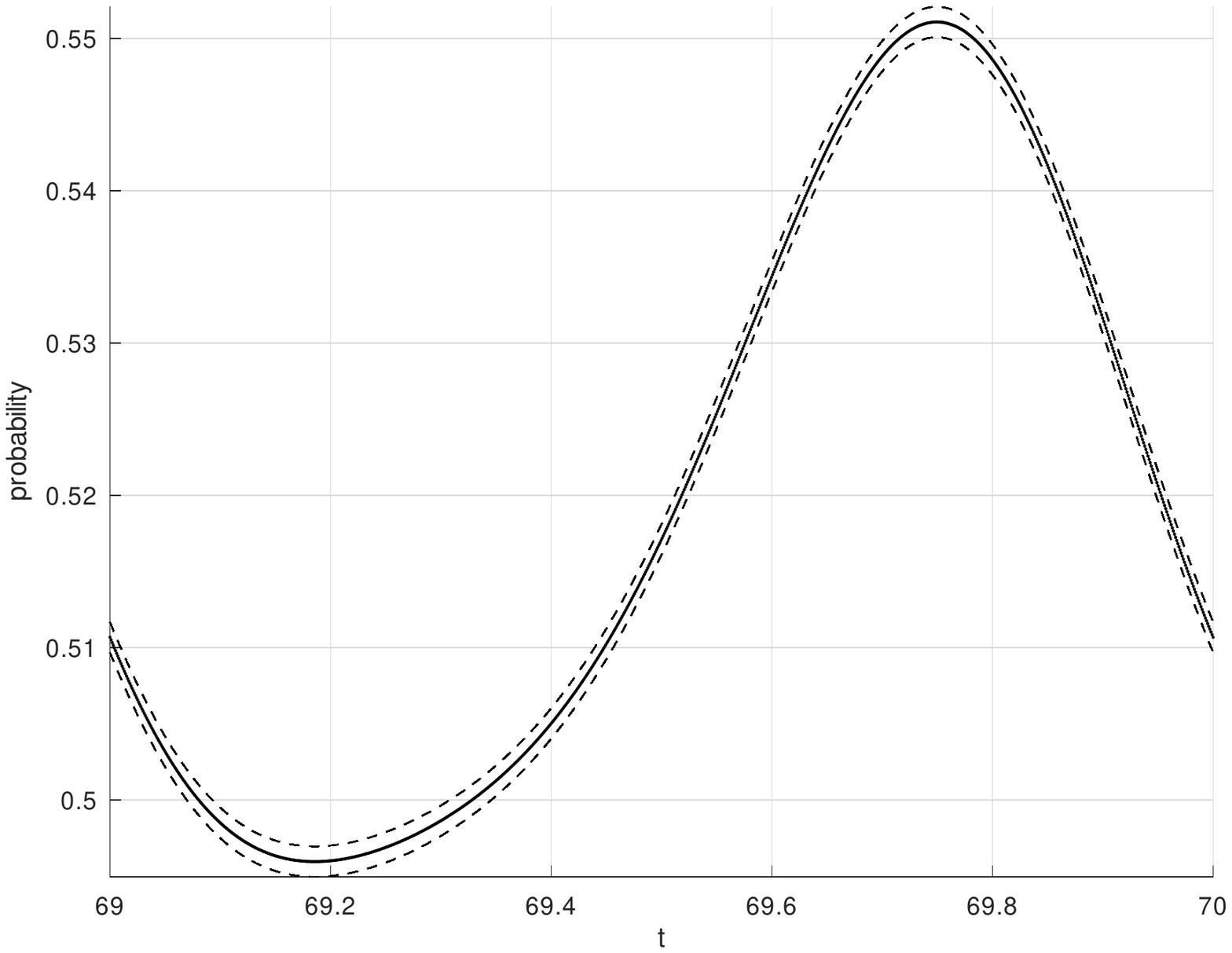}
\vspace{-5cm}\caption{Example 2. Perturbation bounds for the limiting
probability of  empty queue for $ t \in [69,70]$.}
\end{figure}

\begin{figure}
\centering
\includegraphics[height=20cm]{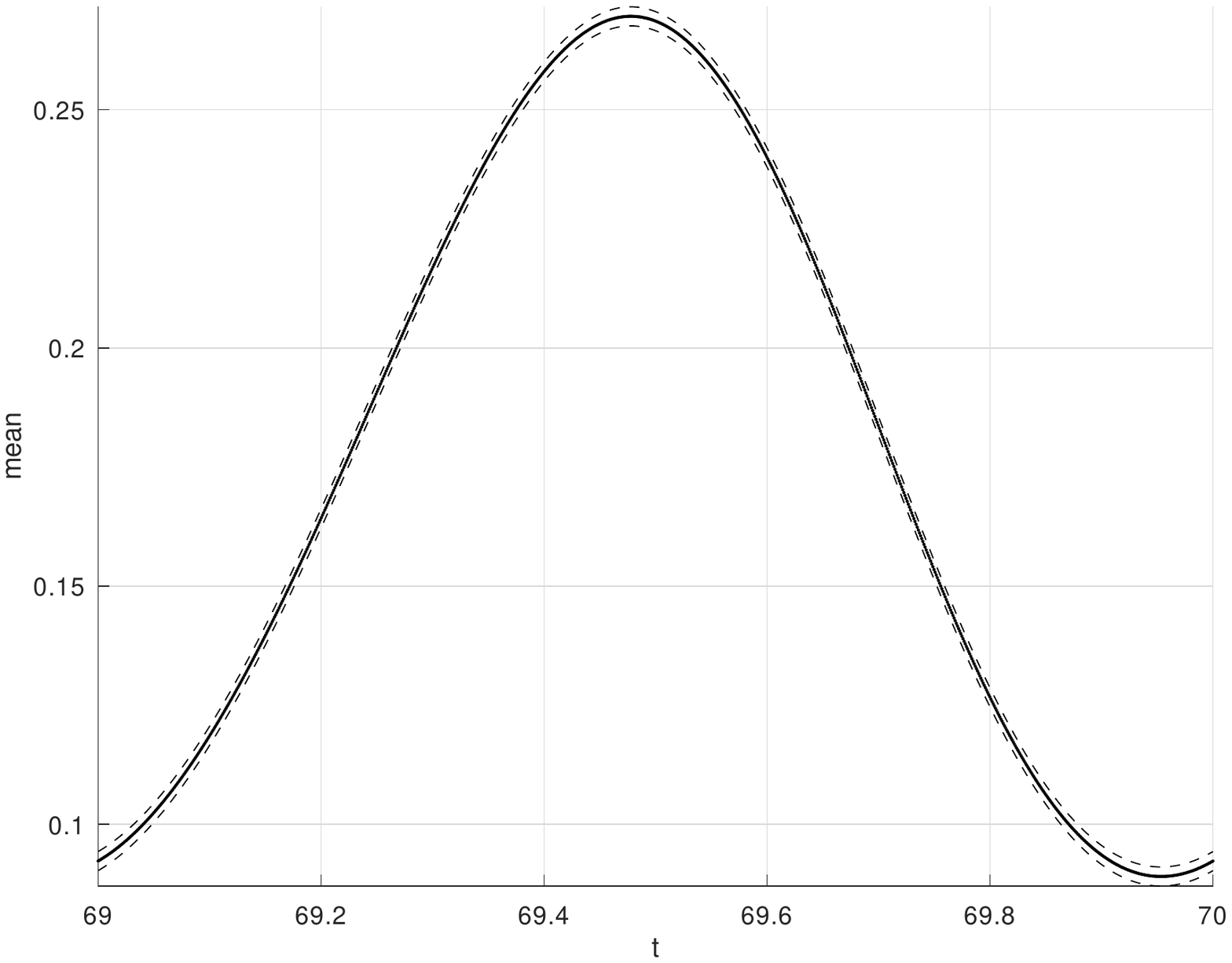}
\vspace{-5cm}\caption{Example 2. Perturbation bounds for  the limiting mean
$E(t,0)$ for $ t \in [69,70]$.}
\end{figure}

\section*{Acknowledgement.} Sections 4 -- 6 were written by A.Z., Y.S. and I.K. under the support of the Russian Science Foundation, project~19-11-00020.

%=============================================================
\end{document}